\documentclass[12pt]{amsart}
\usepackage{amsmath,amsfonts,mathrsfs,amssymb}

\newtheorem{thm}{Theorem}[section]
\newtheorem{theorem}[thm]{Theorem}
\newtheorem{corollary}[thm]{Corollary}
\newtheorem{lemma}[thm]{Lemma}

\theoremstyle{definition}

\theoremstyle{remark}
\newtheorem{remark}[thm]{Remark}

\numberwithin{equation}{section}

\newenvironment{theorem*}[1]{\smallskip\noindent{\bf #1.}\it}{\medskip}

\renewcommand\Im{{\mathrm{Im\,}}}

\newcommand\bN{{\mathbb N}}
\newcommand\bC{{\mathbb C}}
\newcommand\bR{{\mathbb R}}
\newcommand\bZ{{\mathbb Z}}
\newcommand\re{{\mathrm e}}

\newcommand\cF{\mathcal F}

\newcommand\ba{{\mathbf a}}
\newcommand\bb{{\mathbf b}}

\newcommand\be{{\mathbf e}}
\newcommand\bx{{\mathbf x}}
\newcommand\by{{\mathbf y}}
\newcommand\bz{{\mathbf z}}

\newcommand\al{\alpha}
\newcommand\la{\lambda}
\newcommand\eps{\varepsilon}

\begin{document}

\title[Zero asymptotics]{Asymptotics of zeros for some entire
    functions${}^{\dag}$}%
\thanks{${}^{\dag}$The work was partially supported by Ukrainian Foundation
for Basic Research DFFD under  grant No.~01.07/00172.}

\author[R.~O.~Hryniv \and Ya.~V.~Mykytyuk]%
    {Rostyslav O.~Hryniv \and Yaroslav V.~Mykytyuk}%

\address{Institute for Applied Problems of Mechanics and Mathematics,
3b~Naukova st., 79601 Lviv, Ukraine and
Lviv National University, 1 Universytetska st., 79602 Lviv, Ukraine}%
\email{rhryniv@iapmm.lviv.ua, yamykytyuk@yahoo.com}%

\address[Current address of R.H.]{Institut f\"ur Angewandte Mathematik,
Universit\"at Bonn, Wegelerstr. 6, D--53115, Bonn, Germany}
\email{rhryniv@wiener.iam.uni-bonn.de}

\subjclass[2000]{Primary 30D15, Secondary 42A38}%
\keywords{Entire functions, asymptotics of zeros, Fourier transform}%

\date{\today}%
\begin{abstract}
We study the asymptotics of zeros for entire functions of the form
$\sin z + \int_{-1}^1 f(t)\re^{izt}\,dt$ with $f$ belonging to a
space $X \hookrightarrow L_1(-1,1)$ possessing some minimal
regularity properties.
\end{abstract}
\maketitle

\section{Introduction}

The aim of this paper is to study the asymptotics of zeros for
entire functions of the form
\[
    F(z) = m_- \re^{-iz} + m_+ \re^{iz}
        + \int_{-1}^1 f(t)\re^{izt}\,dt,
\]
where $m^\pm$ are nonzero complex numbers and $f$ is an integrable
function over $(-1,1)$. Such a function is the Fourier--Stieltjes
transform~\cite[Ch.~12.5]{Ed} of the measure $m_-\delta_{-1} +
m_+\delta_1 + f(t)\,dt$, with $\delta_x$ being the Dirac measure
at the point~$x$, and, in particular, is an entire function of the
so-called sine type \cite[Lect.~22]{Le}. Asymptotic distribution
of zeros for Fourier--Stieltjes transforms is of key significance
in many areas of function theory, harmonic analysis, functional
analysis etc.~and has been studied for many particular
situations---see, e.g., \cite{KL,Se}. Our interest in this topic
stems, e.g., from the spectral theory of Sturm--Liouville
operators, since the characteristic determinants of such operators
are usually of the above type~(see, e.g., the monographs
\cite{Levit,Ma,PT}).

The standard arguments based on the Riemann--Lebesgue lemma and
Rouch\'e's theorem (see Lemma~\ref{lem:rouche}) show that the
function $F$ has countably many zeros, which can be designated by
$z_n$, $n\in\bZ$, in such a way that $z_n = \pi n + \alpha +
\zeta_n$ with $\alpha:= \frac1{2i}\log(-m_-/m_+)$ and $\zeta_n =
\mathrm{o}(1)$ as $n\to\pm\infty$. We observe that without loss of
generality we may assume that $\alpha=0$; in fact, the function
$\hat F(z):=F(z+\alpha)$ has the form
\[
    \hat F(z) = \hat m_- \re^{-iz} + \hat m_+ \re^{iz}
        + \int_{-1}^1 \hat f(t) \re^{izt}\,dt
\]
with $\hat m_- + \hat m_+=0$ and $\hat f(t):= f(t)\re^{i\alpha
t}$. In other words, there is no loss of generality in putting
$m^\pm:= \pm i/2$, and in this case~$F$ becomes
\begin{equation}\label{eq:Fsin}
    F(z) = \sin z
        + \int_{-1}^1 f(t) \re^{izt}\,dt.
\end{equation}

Since the function $F$ of~\eqref{eq:Fsin} is completely determined
by its zeros $z_n=\pi n + \zeta_n$, a natural problem arises to
describe the sequences $(\zeta_n)_{n\in\bZ}$ that can be obtained
for various functions $f$ from $L_1(-1,1)$ (or from some subset
$X$ in $L_1(-1,1)$). This problem has been attacked in several
particular cases. For instance, the paper by Levin and
Ostrovski\u\i~\cite{LO} contains the following statement.

\begin{theorem*}{Theorem~A}
Assume that $f\in L_p(-1,1)$ with $p\in(1,2]$ and let $z_n = \pi n
+ \zeta_n$ be zeros of the function $F$ in~\eqref{eq:Fsin}. Then
the sequence $(\zeta_n)_{n\in\bZ}$ belongs to $\ell_q$ with
$q:=p/(p-1)$.
\end{theorem*}

It is interesting to know whether the statement of Theorem~A can
be reverted, i.e., whether for any sequence $(c_n)_{n\in\bZ}$ from
$\ell_q$ the numbers $\pi n +c_n$ are zeros (of respective
multiplicities) of the function $F$ in~\eqref{eq:Fsin} with some
$f\in L_p(-1,1)$. Further related questions are about the cases
when $f\in L_p(-1,1)$ with $p=1$ or $p>2$ and, more generally,
when the function $f$ belongs to other spaces embedded into
$L_1(-1,1)$.

It is known that the converse of Theorem~A is true for~$p=2$, cf.
similar statement in~\cite{HMtwo}. It turns out that for
$p\in[1,2)$ this is not true any more: in fact, the numbers
$\zeta_n$ are then Fourier coefficients of some function from
$L_p(-1,1)$ and thus by the Hausdorff--Young
theorem~\cite[Sect.~13.5]{Ed} the corresponding sequences form a
proper subset of $\ell_q$, $q:=p/(p-1)$. Analogous statements
relating smoothness of~$f$ and decay of~$\zeta_n$ were recently
proved in~\cite{HMas} (see also similar results in
\cite{KM,SS})---namely, if $f\in W^s_2(0,1)$, $s\in[0,1]$, then
the $\zeta_n$ are Fourier coefficients of some function from
$W^s_2(0,1)$; cf.~the particular case $s=1$
in~\cite[Lemma~3.4.2]{Ma}. See also~\cite{Zh} for related results
when $f$ is of bounded variation; zero asymptotics of
Fourier--Stieltjes transforms of discrete measures or, more
generally, of almost periodic functions were studied in~\cite{KL}
and \cite[Appendix VI]{Le1}.

The main aim of this note is to point out other function spaces
$X$ that are continuously embedded into $L_1(-1,1)$ and possess
the property that the remainders $\zeta_n$ in the asymptotics of
zeros for the function~$F_f$ are related to Fourier coefficients
of some function~$g\in X$ whenever $f$ belongs to $X$. We show
that such is any Banach space $X\hookrightarrow L_1(-1,1)$ which
satisfies the following assumptions (A1)--(A3):
\begin{itemize}
\item[(A1)] $X$ is an algebra with respect to convolution;
\item[(A2)] the operator $M$ of multiplication by $t$ is
    continuous in~$X$;
\item[(A3)] the set of all finite trigonometric polynomials
    in~$\re^{\pi int}$, $n\in\bZ$, is dense in~$X$.
\end{itemize}
In the following three theorems $X$ denotes a space which is
continuously embedded into $L_1(-1,1)$ and satisfies (A1)--(A3).
We also designate by~$e_n(h)$ the $n$-th Fourier coefficient of a
function $h\in L_1(-1,1)$.

\begin{theorem}\label{thm:main}
Assume that $f\in X$. Then the zeros $z_n=\pi n + \zeta_n$ of the
corresponding entire function $F$ of~\eqref{eq:Fsin} have the
property that $(-1)^n\zeta_n$ are Fourier coefficients of some
function $g$ from~$X$.
\end{theorem}

We remark that if $(-1)^n\zeta_n$ is the $n$-th Fourier
coefficient of a function $g$, then $\zeta_n$ is the $n$-th
Fourier coefficient of the function
\[
    \check{g}(t):=\begin{cases}\ g(t+1),
                \qquad & \text{if}\quad t<0,\\
          \ g(t-1),
                \qquad & \text{if}\quad t\ge0.
    \end{cases}
\]
If therefore the operation $g \mapsto \check{g}$ is one-to-one
on~$X$ (which is, e.g., the case for $L_p(-1,1)$, $p\in
[1,\infty)$, or continuous periodic functions over $[-1,1]$), then
the factor~$(-1)^n$ at $\zeta_n$ in the formulation of
Theorem~\ref{thm:main} can be omitted.

The mapping $f\mapsto g$ established in Theorem~\ref{thm:main} is
multivalued, since different enumerations of zeros $z_n$ generate
different functions~$g$. A natural question arises whether it is
possible to single out single-valued locally continuous branches
of this mapping. We note that such branches do not exist if the
function $F$ of~\eqref{eq:Fsin} has multiple zeros. Indeed, if
$z_0$ is a zero of multiplicity $n\ge2$ and $h\in X$ is such that
$\int_{-1}^1 h(t)\re^{iz_0t}\,dt\ne0$, then the zeros of the
function
    $\sin z  + \int_{-1}^1 [f(t) + w h(t)] \re^{izt}\,dt$
in a neighbourhood of~$z_0$ for sufficiently small complex~$w$ are
given by a Puiseux series in $w$. Therefore there is no
enumeration of these zeros, which is continuous in $w$ close to
$0$ (such an enumeration is only possible on a Riemann surface)
and, as a result, the function $g$ is not then continuous in~$f$.

If, however, all the zeros of the function $F$ in~\eqref{eq:Fsin}
are simple, then such locally continuous branches exist.

\begin{theorem}\label{thm:cont}
Assume that $f_0\in X$ is such that the corresponding function
in~\eqref{eq:Fsin} has only simple zeros. Then there exists a
neighbourhood $\mathcal{O}$ of~$f_0$ and a continuous mapping
$\phi:\mathcal{O} \to X$ such that the numbers  $z_n(f):=\pi n +
(-1)^ne_n(\phi(f))$ are all the zeros of the function $F$
of~\eqref{eq:Fsin} for $f\in\mathcal{O}$.
\end{theorem}

Conversely, we shall show that Fourier coefficients of any $g\in
X$ give zeros of an entire function~\eqref{eq:Fsin} for a unique
$f\in X$.

\begin{theorem}\label{thm:inv}
Assume that $g\in X$. Put $z_n:=\pi n + (-1)^ne_n(g)$, $n\in\bZ$;
then there exists a unique $f\in X$ such that the $z_n$ are zeros
(of respective multiplicities) of the analytic function~$F$
of~\eqref{eq:Fsin}.
\end{theorem}

We observe that for most classical function Banach spaces over
$(-1,1)$ conditions (A1)--(A3) are satisfied and that asymptotics
of Fourier coefficients for functions from such spaces is
relatively well studied, thus providing precise asymptotics of
zeros for the corresponding functions of~\eqref{eq:Fsin}.

For technical reasons it is more convenient for us to make the
change of variables $\phi:\,t\mapsto s:=(t+1)/2$ in the integral
for $F$ and consider the function $\tilde f(s) := 2f(2s-1)$ on
$[0,1]$ instead of the original $f$. The function~$F$ can then be
represented as
\begin{equation}\label{eq:F2t}
   F(z)= F_{\tilde f}(z) := \sin z
        + \int_0^1 \tilde f(s) \re^{iz(2s-1)}\,ds.
\end{equation}
It is clear that the space
  $\tilde X:= \{f \circ \phi^{-1}\mid f\in X\}$
of functions over~$(0,1)$ is a Banach space under the induced norm
$\|\tilde f\|_{\tilde X}:= \tfrac12\|\tilde f \circ \phi\|_X$.
Moreover, $\tilde X$ is continuously imbedded into~$L_1(0,1)$ and
satisfies properties (A1), (A2), while (A3) should be modified to
\begin{itemize}
\item[(A3$'$)] the set of all finite trigonometric polynomials
    in~$\re^{2\pi int}$, $n\in\bZ$, is dense in~$\tilde X$.
\end{itemize}
Observe also that the $n$-th Fourier coefficient of a function
$g\in X$ differs from that of the function $\tilde g:=g \circ
\phi^{-1}$ by the factor $(-1)^n$.

We shall prove Theorems~\ref{thm:main}--\ref{thm:inv} in
Sections~\ref{sec:2} and \ref{sec:3} for the function~$F$ in the
representation~\eqref{eq:F2t} with $\tilde f\in\tilde X$ and with
$(-1)^n e_n(g)$ replaced by $e_n(\tilde g)$. However, to ease the
notation, we shall suppress the tilde and for the rest of the
paper will denote by $X$ a fixed function space over $(0,1)$ that
is continuously imbedded into~$L_1(0,1)$ and satisfies properties
(A1), (A2), and (A3$'$), while $f$ and $g$ will be functions
over~$(0,1)$ from~$X$.

\section{Proofs of Theorems~\ref{thm:main} and
\ref{thm:cont}}\label{sec:2}

Take an arbitrary $f\in X$ and denote by $z_n$ the zeros of the
function $F_f$ of~\eqref{eq:F2t}. The next lemma shows that we can
(and always shall) order $z_n$ in such a way that
    $z_n = \pi n + \mathrm{o}(1)$
as $|n|\to \infty$. Put
\begin{alignat*}{2}
    R_m &:= \{ z\in\bC \mid |z|  \le \pi m + \pi/6\},&
    \qquad &m\in \bN,\\
    K_n &:= \{ z\in\bC \mid |z - \pi n| \le \pi/6\},&
    \qquad &n\in\bZ.
\end{alignat*}

\begin{lemma}\label{lem:rouche}
Assume that $\mathcal{O}$ is a neighbourhood of $f_0\in X$ such
that
\[
    \mathcal O \subset \{f\in X \mid \|f-f_0\|_{L_1}\le 1/8\}.
\]
Then there exists $n_0\in \bN$ such that for every $f\in\mathcal
O$ the function $F_f$ has precisely $2m+1$ zeros in the disk $R_m$
for $m>n_0$ and precisely one zero in $K_n$ for $|n|>n_0$.
\end{lemma}

\begin{proof}
We recall that
\begin{equation}\label{eq:rouche}
    \int_{0}^1  f(t) \re^{iz(2t-1)}\,dt = \mathrm{o}(\re^{|\Im z|})
\end{equation}
as $|z|\to\infty$  for all $f\in L_1(0,1)$ \cite[Lemma~1.3.1]{Ma}.
Denote by $\partial R_m$ the boundary of $R_m$; it follows then
from~\eqref{eq:rouche} that there exists $n_0\in\bN$ such that,
for all $m>n_0$,
\[
    \max_{z\in \partial R_m}
    \Bigl| \re^{-|\Im z|}\int_{0}^1  f_0(t) \re^{iz(2t-1)}\,dt\Bigr|
    < \frac18,
\]
so that
\[
    \max_{z\in \partial R_m}
    \Bigl| \re^{-|\Im z|}\int_{0}^1  f(t) \re^{iz(2t-1)}\,dt\Bigr|
    < \frac14
\]
for all $f\in \mathcal O$. It is readily seen that
\begin{equation}\label{eq:2.Rm}
    \min_{z\in \partial R_m} \bigl|\re^{-|\Im z|}\sin z \bigr| \ge
    \frac14
\end{equation}
for all $m\in \bN$. In fact, writing $z= x + iy$ with $x,y\in\bR$
and taking $y\ge 0$ for the sake of definiteness, we get
\[
    \bigl|\re^{-|\Im z|}\sin z\bigr|
        = \Bigl|\re^{-y}\frac{\re^{ix-y}- \re^{-ix+y}}{2i}\Bigr|
        = \Bigl|\frac{\re^{-2y}- \re^{-2ix}}{2}\Bigr|.
\]
Since $\bigl|\re^{-2ix}\bigr|=1$, we get
 $|\re^{-2y}- \re^{-2ix}|\ge\tfrac12$ as soon as $y\ge\tfrac12\log2$. For
$y<\tfrac12\log2$ and $z\in\partial R_m$ one concludes that
 \(
    |x| = \sqrt{(\pi m + \pi / 6)^2 - y^2}
        > \pi m + \pi/12,
 \)
hence $|2x| \in (2\pi m + \pi /6, 2\pi m + \pi/3)$ and \(
\bigl|\Im \re^{-2ix}\bigr|= |\sin (2x)| \ge \tfrac12 \), so that
$|\re^{-2y}- \re^{-2ix}|\ge\tfrac12$. Thus~\eqref{eq:2.Rm} holds,
and Rouch\'e's theorem yields now the conclusion about~$R_m$.
Inclusion of zeros in $K_n$ is justified in the same manner.
\end{proof}

 Put $\tilde z_n := z_n - \pi n$; then the relation
\[
    F_f(\pi n + \zeta) =
       (-1)^n \sin \zeta + (-1)^n
       \int_0^1 f(t) \re^{i\zeta(2t-1)}\re^{-2\pi nit}\,dt
\]
shows that the numbers $\tilde z_n$ satisfy the equality
\[
    \sin \tilde z_n +
       \int_0^1 f(t) \re^{i\tilde z_n(2t-1)}\re^{-2\pi nit}\,dt =0.
\]
Developing $\re^{i\tilde z_n(2t-1)}$ into the Taylor series and
then changing the summation and integration order, we get the
following equation for~$\tilde z_n$:
\[
     \sin\tilde z_n
         +\sum_{k=0}^\infty\frac{\tilde z_n^{k}}{k!}
            \int_0^1 i^k(2t-1)^{k} f(t) \re^{-2\pi nit}\,dt=0.
\]
Denote by $M$ the operator of multiplication by $i(2t-1)$; then
the above equation can be recast as
\begin{equation}\label{eq:zeta}
    \sin\tilde z_n  +  \sum_{k=0}^\infty
            \frac{e_n(M^{k}f)\tilde z_n^{k}}{k!} =0,
\end{equation}
where
\[
    e_n(h):= \int_0^1 h(t)\re^{-2\pi nit}\, dt
\]
is the $n$-th Fourier coefficient of a function $h\in L_1(0,1)$.

Our next goal is to show that the numbers $\tilde z_n$ form a
sequence of Fourier coefficients of some function~$g$ from~$X$. To
this end we shall treat the system of equations~\eqref{eq:zeta}
indexed by $n\in\bZ$ as a single equation in some Banach space for
the sequence $\bz=(\tilde z_n)_{n\in\bZ}$.

We denote by $\cF$ the discrete Fourier transform $\cF: L_1(0,1)
\to \ell_\infty(\bZ)$ given by
\[
    \cF f:= \be(f):= \bigl( e_n(f)\bigr)_{n\in\bZ},
\]
and write $\widehat X$ for the linear manifold $\cF X \subset
\ell_\infty$. Being endowed with the norm \[
        \|\ba\|_{\widehat X} := \|\cF^{-1}\ba\|_X,
\]
the set $\widehat X$ becomes a Banach space. For elements
$\ba=(a_n)$ and $\bb=(b_n)$ of $\widehat X$, we denote by $\ba\bb$
the entrywise product of $\ba$ and $\bb$, i.e., the element of
$\ell_\infty$ with the $n$-th entry $a_nb_n$. The next lemma shows
that this multiplication is continuous in~$\widehat X$.

\begin{lemma}
There exists $\rho>0$ such that $\|\ba\bb\|_{\widehat X}\le\rho
\|\ba\|_{\widehat X}\|\bb\|_{\widehat X}$.
\end{lemma}

\begin{proof}
We observe that $\cF^{-1} (\ba\bb)=
(\cF^{-1}\ba)\ast(\cF^{-1}\bb)$, where $f\ast g$ stands for the
convolution of $L_1$-functions $f$ and $g$,
\[
    (f\ast g)(x) = \int_0^1 f(x-t) g(t)\,dt.
\]
Here we assume that the function $f$ is periodically extended to
the interval $(-1,0)$. It remains to recall that by (A1) the
convolution operation is continuous in $X$.
\end{proof}

\begin{corollary}\label{cor:powers} If $\ba\in\widehat X$, then
$\ba^n\in\widehat X$ for all $n\in\bN$ and, moreover, with the constant
$\rho$ of the previous lemma, one has $\|\ba^n\|_{\widehat X}\le
\rho^{n-1}\|\ba\|^n_{\widehat X}$.
\end{corollary}

\begin{corollary}\label{cor:contin} Suppose that $\ba,
\bb\in\widehat X$ and put $a:=\max\{\|\ba\|_{\widehat X},\|\bb\|_{\widehat X}\}$.
Then, for all $n\in\bN$,
\[
    \|\ba^{n} - \bb^{n}\|_{\widehat X}
            \le n (\rho a)^{n-1} \|\ba-\bb\|_{\widehat X}.
\]
\end{corollary}

We denote by $\Gamma$ the linear manifold of sequences
$\gamma :=(\ba_k)_{k=0}^\infty$ in $\widehat X$, for which the following norm
$\|(\ba_k)\|_\Gamma$ is finite:
\[
    \|(\ba_k)\|_\Gamma := \|\ba_0\|_{\widehat X} +
           \sum_{k=1}^\infty \frac{\|\ba_{k}\|_{\widehat X}}{(k-1)!}.
\]
It is easily verified that $\Gamma$ becomes a Banach space under the
above norm. With every $\gamma =(\ba_k)\in\Gamma$ we associate a nonlinear
mapping $G_\gamma \,:\,\widehat X \to \widehat X$ given by
\[
    G_\gamma (\bx) := \sum_{k=1}^\infty\frac{(-1)^{k}\bx^{2k+1}}{(2k+1)!}
            - \ba_0 - \sum_{k=1}^\infty\frac{\ba_k\bx^k}{k!}.
\]
Observe that the system of equations~\eqref{eq:zeta} assumes the
form $\bx= G_\gamma (\bx)$ for a particular choice of the elements
$\ba_k$, namely for $\ba_{k}=\be(M^{k}f)$. We shall show that for
every $\gamma \in\Gamma$ of sufficiently small norm the mapping
$G_\gamma $ is contractive in some ball $B$ of $\widehat X$
centered at the origin and thus the equation $\bx = G_\gamma
(\bx)$ has a unique solution in $B$. In the following, $B_Y(r)$
will stand for the closed ball of a Banach space~$Y$ with center
at the origin and radius~$r$.

\begin{lemma}\label{lem:contr}
Suppose that $\gamma \in\Gamma$ and $r>0$ are chosen so that
$2\|\gamma \|_\Gamma\le r\le(2\rho)^{-1}$. Then the mapping $G_\gamma $ leaves the
ball $B_{\widehat X}(r)$ invariant and for every $\bx,\by\in B_{\widehat X}(r)$
we have
\[
    \|G_\gamma (\bx) - G_\gamma (\by)\|_{\widehat X}
        \le \tfrac12 \|\bx-\by\|_{\widehat X}.
\]
\end{lemma}

\begin{proof}
Recalling the definition of the mapping $G_\gamma $ and using
Corollaries~\ref{cor:powers} and \ref{cor:contin} and the
inequality $|\sinh x  - x|\le x/2$ for $x\in[0,\tfrac12]$, we find
that, for any $\bx\in B_{\widehat X}(r)$ with $r$ as in the
statement of the lemma,
\begin{align*}
    \|G_\gamma  (\bx)\|_{\widehat X} &\le
        \sum_{k=1}^\infty\frac{\rho^{2k}r^{2k+1}}{(2k+1)!}
    + \sum_{k=0}^\infty
            \frac{(\rho r)^{k}\|\ba_{k}\|_{\widehat X}}{k!}\\
    &\le \frac1\rho\bigl[\sinh(\rho r) - \rho r\bigr] + \|\gamma \|_{\Gamma}
        \le \tfrac{r}2  +  \|\gamma \|_{\Gamma} \le r,
\end{align*}
i.e., that $G_\gamma  B_{\widehat X}(r) \subset B_{\widehat X}(r)$. In a similar manner
we estimate the $\widehat X$-norm of the difference $G_\gamma (\bx)-G_\gamma (\by)$ as
\begin{align*}
    \|G_\gamma (\bx) - G_\gamma (\by)\|_{\widehat X}
    &\le \sum_{k=1}^\infty
        \frac{\|\bx^{2k+1}-\by^{2k+1}\|_{\widehat X}}{(2k+1)!}
    +\sum_{k=1}^\infty
        \frac{\|\ba_{k}(\bx^{k}-\by^{k})\|_{\widehat X}}{k!}\\
       & \le  \sum_{k=1}^\infty \frac{(\rho r)^{2k}
                \|\bx-\by\|_{\widehat X}}{(2k)!}
        + \rho\sum_{k=1}^\infty \frac{(\rho r)^{k-1}
         \|\bx-\by\|_{\widehat X}\|\ba_{k}\|_{\widehat X}}{(k-1)!}\\
       & \le \bigl[\cosh (\rho r) -1 + \rho \|\gamma \|_\Gamma \bigr]
           \|\bx-\by\|_{\widehat X}\\
       & \le \rho\bigl(\tfrac{r}2 + \|\gamma \|_\Gamma \bigr)
           \|\bx-\by\|_{\widehat X}
        \le \tfrac12\|\bx-\by\|_{\widehat X},
\end{align*}
and the proof is complete.
\end{proof}

The Banach fixed point principle yields now the following

\begin{corollary}\label{cor:contr}
Let $r_0:=(4\rho)^{-1}$. Then for every $\gamma $
in the ball $B_\Gamma(r_0)$ the equation $\bx= G_\gamma (\bx)$ has a unique
solution $\bx=\bx(\gamma )$ in the ball $B_{\widehat X}(2r_0)$.
\end{corollary}

We shall show next that the solution $\bx(\gamma )$ given in the above
corollary depends continuously on $\gamma $.

\begin{lemma}\label{lem:cont}
With $r_0:=(4\rho)^{-1}$, the map $B_\Gamma(r_0) \ni \gamma
\mapsto \bx(\gamma )\in B_{\widehat X}(2r_0)$ is continuous.
\end{lemma}

\begin{proof} Assume that $\gamma=(\ba_k)$ and
$\tilde \gamma= (\tilde\ba_k)$ belong to $B_\Gamma(r_0)$; then for
an arbitrary $\bx\in B_{\widehat X}(2r_0)$ the following
inequality holds:
\[
    \|G_\gamma (\bx) - G_{\tilde \gamma }(\bx)\|_{\widehat X}
     \le \|\ba_0-\tilde\ba_0\|_{\widehat X} +
         \sum_{k=1}^\infty \frac{\rho^k\|\bx\|^k_{\widehat X}
         \|\ba_k-\tilde\ba_k\|_{\widehat X}}{k!}
     \le \|\gamma -\tilde \gamma \|_\Gamma.
\]
Bearing in mind the relations
\[
    \bx(\gamma ) = G_\gamma \bigl(\bx(\gamma )\bigr), \qquad
    \bx(\tilde \gamma ) = G_{\tilde \gamma }\bigl(\bx(\tilde \gamma )\bigr),
\]
we find that
\begin{align*}
    \|\bx(\gamma ) - \bx(\tilde \gamma )\|_{\widehat X} &\le
        \|G_\gamma \bigl(\bx(\gamma )\bigr)-
            G_{\tilde \gamma }\bigl(\bx(\gamma )\bigr)\|_{\widehat X}
    +   \|G_{\tilde \gamma }\bigl(\bx(\gamma )\bigr)-
            G_{\tilde \gamma }\bigl(\bx(\tilde \gamma )\bigr)\|_{\widehat X} \\
    & \le \|\gamma -\tilde \gamma \|_\Gamma +
        \tfrac12\|\bx(\gamma ) - \bx(\tilde \gamma )\|_{\widehat X},
\end{align*}
so that $\|\bx(\gamma) - \bx(\tilde\gamma)\|_{\widehat X} \le
2\|\gamma -\tilde \gamma \|_\Gamma$, and the proof is complete.
\end{proof}

\begin{proof}[Proof of Theorem~\ref{thm:main}]
As we have observed, the system of equations~\eqref{eq:zeta}
assumes the form $\bx=G_\gamma (\bx)$ for the element $\gamma
=(\ba_k)_{k\in\bZ_+}$ of $\Gamma$ with $\ba_k:=\be(M^k f)$. Since
$\|\gamma \|_\Gamma \le \bigl[1+ \exp(\|M\|)\bigr]\|f\|_X$, we
conclude by Corollary~\ref{cor:contr} that the equation $\bx=
G_\gamma (\bx)$ has a solution $\bx(\gamma)\in B_{\widehat
X}(2r_0)$ if $\|f\|_X$ is small enough. In a general situation we
shall replace $\gamma $ with $\tilde \gamma \in\Gamma$ in such a
manner that $\tilde \gamma $ is of sufficiently small norm and the
solution $\tilde \bx= (\tilde x_n)_{n\in\bZ}$ of the equation
$\bx=G_{\tilde \gamma }(\bx)$ is such that $\tilde x_n = \tilde
z_n$ if $|n|$ is large enough. Modifying the function
$\cF^{-1}\tilde \bx\in X$ by a trigonometric polynomial, we get
then a function $g\in X$ such that $e_n(g) = \tilde z_n$ for all
$n\in\bZ$.

The details are as follows. As $M$ is bounded in $X$, there exists
$k_0\in\bN$ such that
\[
    \sum_{k>k_0}\frac{\|\be(M^kf)\|_{\widehat X}}{(k-1)!}
         < \frac1{20\rho}.
\]
Since by assumption (A3$'$) the trigonometric polynomials in
$\re^{2\pi inx}$ form a dense set in $X$, for each $k=0,\dots k_0$
there exists a polynomial $p_k$ of degree $d_k$ such that
$\|M^kf-p_k\|_X<(20\rho)^{-1}$. We set now $\tilde \gamma  =
(\tilde \ba_k)_{k\in\bZ_+}$ with $\tilde \ba_k := \be(M^kf-p_k)$
for $k=0,\dots,k_0$, and $\tilde \ba_k := \be(M^kf)$ for $k>k_0$.
Then $\|\tilde \gamma \|_\Gamma < r_0=(4\rho)^{-1}$, so that by
Corollary~\ref{cor:contr} the equation $\bx=G_{\tilde \gamma
}(\bx)$ has a unique solution $\tilde \bx= (\tilde x_n)_{n\in\bZ}$
in $B_{\widehat X}(2r_0)$.

We observe that $e_n(p_k)=0$ if $k=0,\dots,k_0$ and
$|n|>d:=\max\{d_0,\dots,d_{k_0}\}$. Thus the $n$-th component of
the relation $\tilde \bx=G_{\tilde \gamma }(\tilde\bx)$ for such
$n$ reads
\[
    \sin\tilde x_n  +  \sum_{k=0}^\infty
            \frac{e_n(M^{k}f)\tilde x_n^{k}}{k!} =0.
\]
Lemma~\ref{lem:rouche} gives $n_0\in\bN$ such that
equation~\eqref{eq:zeta} has a unique solution in the
$\pi/6$-neighbourhood of the origin for all $n\in\bZ$ with
$|n|>n_0$. Since $\tilde \bx\in\widehat X$ implies
$\lim_{|n|\to\infty}\tilde x_n=0$, there exists a number $n_1>d$
such that $\tilde z_n=\tilde x_n$ if $|n|>n_1$. Then the function
\begin{equation}\label{eq:g}
    g :=  \cF^{-1}\tilde \bx + \sum_{n=-n_1}^{n_1}(\tilde z_n -\tilde x_n)
        \re^{2\pi inx}
\end{equation}
belongs to $X$ and satisfies the equalities $e_n(g)=\tilde z_n$
for all $n\in\bZ$, and the proof is complete.
\end{proof}

\begin{proof}[Proof of Theorem~\ref{thm:cont}]
Assume that $f_0$ is such that all zeros of $F_{f_0}$ are simple.
The implicit function theorem then states that for every
zero~$\zeta_0$ of $F_{f_0}$ there exists a neighbourhood $\mathcal
{O}_{\zeta_0}$ of $f_0$ and a continuous function $\zeta:\mathcal
O_{\zeta_0} \to \bC$ such that $\zeta(f)$ is a zero of $F_f$ for
all $f\in\mathcal O_{\zeta_0}$ and $\zeta(f_0)=\zeta_0$. In view
of formula~\eqref{eq:g} it suffices to show that there exists a
neighbourhood $\mathcal {O}$ of $f_0$ such that the element
$\tilde \bx=\tilde \bx(f)\in\widehat X$ can be taken continuous in
$f\in\mathcal O$ and the number $n_1$ in~\eqref{eq:g} can be taken
independently of $f\in\mathcal O$.

Construct the element $\tilde \gamma (f_0)=
(\ba_k(f_0))_{k\in\bZ_+}$ as in the proof of
Theorem~\ref{thm:main}; in particular, there is a number $k_0\in
\bN$ and polynomials $p_k$ of degree $d_k$, $k=0,\dots ,k_0$, such
that $\tilde \ba_k(f_0) = \be(M^kf_0-p_k)$ for $k=0,\dots,k_0$,
and $\tilde \ba_k(f_0) = \be(M^kf_0)$ for $k>k_0$. We put now
$\tilde \gamma (f)=(\ba_k(f))_{k\in\bZ_+}$ with
$\ba_k(f)=\ba_k(f_0) + \be\bigr(M^k(f-f_0)\bigr)$; then
\[
    \|\tilde \gamma (f) - \tilde \gamma (f_0)\|_\Gamma
        \le (1+\exp(\|M\|))\|f-f_0\|_{X}
\]
and thus $\tilde \gamma (f)$ is continuous in $f$. In particular,
$\|\tilde \gamma (f)\|_\Gamma \le (4\rho)^{-1}$ if $\|f-f_0\|_X$
is smaller than a certain $\eps_0>0$, so that by
Corollary~\ref{cor:contr} the equation $\bx=G_{\tilde \gamma
(f)}(\bx)$ has then a unique solution $\tilde \bx(f)= (\tilde
x_n(f))_{n\in\bZ}$ in $\widehat X$. By Lemma~\ref{lem:cont} this
solution $\tilde \bx(f)$ depends continuously on $\tilde \gamma
(f)$ and thus on $f\in X$ with $\|f-f_0\|_X<\eps_0$.

Now we see that the numbers $\tilde x_n(f)$ for
$|n|>d=\max\{d_0,\dots,d_{k_0}\}$ satisfy the equation
\[
    \sin\tilde x_n  +  \sum_{k=0}^\infty
            \frac{e_n(M^{k}f)\tilde x_n^{k}}{k!} =0.
\]
It follows from Lemma~\ref{lem:rouche} that $\pi n + \tilde
x_n(f)$ is a zero of $F_f$ as soon as $|\tilde x_n(f)|<\pi/6$ and
$|n|>n_0$ with $n_0$ of that lemma. Observe that the inclusion
$X\hookrightarrow L_1(0,1)$ implies that there exists a constant
$C>0$ such that $\|h\|_{L_1}\le C \|h\|_X$ for all $f\in X$ and
hence for all $\bx = (x_n)$ and $\by = (y_n)$ in $\widehat X$ and
all $n\in\bZ$ it holds
\[
    |x_n-y_n|=|e_n(\cF^{-1}(\bx-\by))|\le \|\cF^{-1}(\bx-\by)\|_{L_1}
        \le C\|\bx-\by\|_{\widehat X}.
\]
In particular, $\tilde x_n(f)$ is a continuous function of~$f$ if
$\|f-f_0\|_X<\eps_0$. Take now $n_0'\in\bN$ such that $|\tilde
x_n(f_0)|<\pi/12$ for all $n\in\bZ$ with $|n|>n_0'$. Then, by the
above remark, $|\tilde x_n(f)|<\pi/6$ for all $n\in\bZ$ with
$|n|>n_0'$ if $\|f-f_0\|_X$ is smaller than a certain
$\eps_1\in(0,\eps_0)$. Taking $\mathcal O'$ to be the
$\eps_1$-neighbourhood of $f_0$ in $X$ and
$n_1:=\max\{d,n_0,n_0'\}$, we conclude that $\pi n +\tilde x_n(f)$
is a zero of $F_f$ for all $f\in\mathcal O'$ and all $n\in\bZ$
with $|n|>n_1$.

We enumerate now the zeros $z_n(f)$ of $F_{f}$ for $f\in\mathcal
O$ in such a manner that $\tilde z_n(f):=z_n(f)-\pi n$ coincides
with $\tilde x_n(f)$ if $|n|>n_1$. By the remark made at the
beginning of the proof there exists a neighbourhood $\mathcal O''$
of $f_0$ such that the zeros $z_n(f)$ can be chosen continuously
for $n=-n_1,\dots,n_1$ if $f$ varies over $\mathcal O''$. It
remains to take
\[
    \phi(f) := \cF^{-1} \tilde\bx(f) +
        \sum_{n=-n_1}^{n_1}(\tilde z_n(f) -\tilde x_n(f))
        \re^{2\pi inx}
\]
for $f\in\mathcal O :=\mathcal O'\cap \mathcal O''$, and the proof
is complete.
\end{proof}

\section{Proof of Theorem~\ref{thm:inv}}\label{sec:3}

Assume that $g\in X$, and denote by
$A_g$ a linear operator in $X$ given by
\[
    A_g (f) = f + \cF^{-1}\sum_{k=1}^\infty
        \frac{\be(M^kf)\be^k(g)}{k!},
\]
where, as above, $M$ stands for the operator of multiplication by
$i(2t-1)$.

\begin{lemma}\label{lem:4.1}
The operator $A_g$ is bounded in $X$; moreover, it is boundedly
invertible if $\|g\|_X\le (2\rho \|M\|)^{-1}$.
\end{lemma}

\begin{proof}
In view of Corollaries~\ref{cor:powers} and \ref{cor:contin}, we find
that
\[
 \|(A_g-I)f\|_X \le \sum_{k=1}^\infty
        \frac{\rho^k\|M\|^k\|g\|^k_X\|f\|_X}{k!}
        = \bigl[\exp\bigl(\rho\|M\|\|g\|_X\bigr)-1 \bigr]\|f\|_X,
\]
so that $A_g$ is bounded. If, moreover, $\|g\|_X\le (2\rho
\|M\|)^{-1}$, then $\|A_g-I\|<1$ and hence $A_g$ is boundedly
invertible in $X$.
\end{proof}

\begin{corollary}\label{cor:4.2}
Assume that  $g\in X$ is such that $\|g\|_X\le
(2\rho \|M\|)^{-1}$ and $\|g\|_{L_1}<\pi/2$. Put $z_n:=\pi n +
e_n(g)$, $n\in\bZ$; then there exists a unique function $f\in X$
such that $F_f(z_n) = 0$ for all $n\in\bZ$.
\end{corollary}

\begin{proof}
We put
\[
    \ba:= \sin \be(g) :=  \sum_{k=0}^\infty
        (-1)^{k}\frac{\be^{2k+1}(g)}{(2k+1)!};
\]
it is clear that the series converges in $\widehat X$, so that
$\ba\in\widehat X$. By Lemma~\ref{lem:4.1} the operator~$A_g$ is
invertible in $X$, so that the function
 $f:= - A_g^{-1}(\cF^{-1}\ba)$
is well defined and belongs to $X$.

The definition of the operator~$A_g$ implies that for each
$n\in\bN$ it holds
\[
    \sin e_n(g) + \sum_{k=0}^\infty
        \frac{e_n(M^kf)e_n^k(g)}{k!} =0.
\]
Reverting the arguments of the previous section, we see that the
function $F_f$ vanishes at the points $z_n$ (which are pairwise
distinct).

To prove uniqueness, assume that there is $\tilde f\in X$ such
that $F_{\tilde f}(z_n)=0$. We put $h:=f-\tilde f$ and
$H(z):=\int_0^1 h(t) \re^{iz(2t-1)}\,dt$ and observe that
$H(z_n)=0$ for all $n\in\bZ$. Writing $H(z)$ as $\int_{-1}^1 q(t)
\re^{izt}\,dt$ with $q(t):=\frac12h(\frac{1+t}2)$ and applying
Lemma~\ref{lem:3.3}, we conclude that $h = 0$, i.e., that $f =
\tilde f$. The corollary is proved.
\end{proof}

The following uniqueness statement follows from Levinson's
theorems on closure~$L_p(-1,1)$ of a system of
exponentials~\cite[Chap.~I]{Levins} for the case $p>1$; however,
these theorems are not applicable in our case $p=1$.

\begin{lemma}\label{lem:3.3}
Assume that $q\in L_1(-1,1)$ and that the function
\[
    Q(z):= \int_{-1}^1 q(t)\re^{izt}\,dt
\]
has zeros $z_k$, $k\in\bZ$ (of respective multiplicity) of the
form $z_k= \pi k + \tilde z_k$, where $\tilde z_k$ are Fourier
coefficients of a function $r$ from $L_1(-1,1)$. Then $q=0$.
\end{lemma}

\begin{proof}
Without loss of generality we may assume that none of $z_k$
vanishes. Indeed, by~\cite[Theorem~VI]{Levins} finitely many of
$z_k$ may be arbitrarily changed without affecting the property of
the set $\{\re^{iz_k t}\}_{k\in\bZ}$ to be closed~$L_1(-1,1)$.
(Here, as usual, if some $\omega$ is repeated $l>1$ times among
$z_k$, then we augment the above system by the functions
$t\re^{i\omega t}, \dots, t^{l-1}\re^{i\omega t}$.)

Assume that $q$ is not identically zero. Observe that the numbers
$z_n$ give then all the zeros of the entire function~$Q$. Indeed,
simple arguments based on the Jensen's theorem lead to the
inequality
\begin{equation}\label{eq:3.jen}
    \int_0^r \frac{n(t)}t\,dt \le \frac{2r}\pi + C,
\end{equation}
where $n(r)$ is the number of zeros of~$Q$ in the disc
$\{z\in\bC\mid |z|\le r\}$ and $C$ is some constant. If, however,
$Q$ had at least one additional zero, then for any $\eps>0$ and
all sufficiently large $t>0$ we would have the estimate
\[
    n(t) \ge 2 \Bigl[\frac{t-\eps}\pi\Bigr] + 2
\]
($[\,s\,]$ denoting the integral part of a number $s\in\bR$),
which contradicts~\eqref{eq:3.jen}. Alternatively, we could then
apply the analogue of Theorem~III from~\cite{Levins} for
exponentials $\re^{i\la_k t}$ with complex $\la_k$ to contradict
the assumption that $q$ is not identically zero.

Since the function $Q$ is clearly of exponential type not
exceeding $1$, it can be represented by its Hadamard canonical
product as
\[
    Q(z) = \re^{az+b} \lim_{n\to\infty} \prod_{k=-n}^n
        \Bigl( 1 - \frac{z}{z_k}\Bigr)\re^{z/z_k}
\]
with some constants $a$ and $b$. Observe that the series
$\mathrm{V.p.}\!\sum_{k=-\infty}^\infty 1/{z_k}$ converges, hence
the factors $\re^{z/z_k}$ can be incorporated into $\re^{az+b}$ by
modifying~$a$ accordingly. Writing $\sin z$ as the canonical
product, we conclude that, for $z\ne \pi k$, $k\in\bZ$,
\[
    \frac{Q(z)\re^{-az-b}}{\sin z}
        = \Bigl(\frac1z-\frac1{z_0}\Bigr)
            \mathrm{V.p.}\! \prod_{k=-\infty}^\infty\!\!'
         \Bigl( 1 - \frac{\tilde z_k}{z_k}\Bigr)
             \mathrm{V.p.}\!\prod_{k=-\infty}^\infty\!\!'
         \Bigl( 1 - \frac{\tilde z_k}{z-\pi k}\Bigr),
\]
where the prime means that the factor corresponding to $k=0$ is
omitted. Both products above converge due to the convergence of
the series
 \(
    \mathrm{V.p.}\!\sum_{-\infty}^{'\infty} \tilde z_k/k,
 \)
see~\cite[10.1.5]{Ed}.

We denote
\begin{align*}
    Q_1(iy) &:= \mathrm{V.p.}\!\prod_{k=-\infty}^\infty\!\!'
         \Bigl( 1 - \frac{\tilde z_k}{iy-\pi k}\Bigr)\\
        & = \prod_{k=1}^{\infty} \Bigl(
            1 + \frac{iy(\tilde z_k + \tilde z_{-k})}{y^2 + \pi^2k^2}
              + \frac{\pi k(\tilde z_k - \tilde z_{-k})}{y^2 + \pi^2k^2}
              - \frac{\tilde z_k \tilde z_{-k}}{y^2 + \pi^2k^2}
              \Bigr)
\end{align*}
and prove that $\lim_{y\to\pm\infty}Q_1(iy)=1$. To this end it
suffices to show that the sums of the three series
\[
    \sum_{k=1}^\infty \frac{y(\tilde z_k + \tilde z_{-k})}{y^2 +
    \pi^2k^2},
        \qquad
    \sum_{k=1}^\infty \frac{\pi k(\tilde z_k - \tilde z_{-k})}
        {y^2 + \pi^2k^2},       \qquad
    \sum_{k=1}^\infty \frac{\tilde z_k \tilde z_{-k}}{y^2 + \pi^2k^2}
\]
vanish as $y\in\bR$ tends to $\pm\infty$. For the second and third
series this follows from their uniform in $y\in\bR$ convergence.
Indeed, the second series converges uniformly on $\bR$ by the
Abel--Dirichlet test (to this end we rewrite the $k$-th summand as
    $\frac{\tilde z_k - \tilde z_{-k}}{\pi k}
     \frac{\pi^2k^2}{y^2 + \pi^2k^2}$
and recall that the series
    $\sum_{k=1}^\infty(\tilde z_k - \tilde z_{-k})/k$
converges~\cite[10.1.5]{Ed}), while to the third one the dominated
convergence test applies. Observe that $y/(y^2+\pi^2k^2)$ for
$y\ne0$ is the $k$-th Fourier coefficient of the
function~$u_y(t):= \cosh(y-y|t|)/(2\sinh y)$. Therefore $y\tilde
z_k/(y^2+\pi^2k^2)$ is the $k$-th Fourier coefficient of the
convolution $w_y:= u_y \ast r$. The function $w_y$ is in
$W_2^1(-1,1)$, hence its Fourier series converges pointwise to
$w_y$; in particular, we have
\begin{equation}\label{eq:3.series}
    \mathrm{V.p.}\frac12\sum_{k=-\infty}^\infty
        \frac{y \tilde z_k}{y^2 + \pi^2k^2} = w_y(0)
        = \int_{-1}^1 u_y(-t) r(t)\,dt.
\end{equation}
Since the functions $u_y$ are bounded on $(-1,1)$ uniformly in
$y\in\bR\setminus(-1,1)$ and go to zero pointwise for all nonzero
$t\in(-1,1)$ as $|y|\to\infty$, the Lebesgue dominated convergence
theorem shows that the integral in~\eqref{eq:3.series} vanishes as
$|y|\to\infty$.

The above considerations show that
\[
    \lim_{y\to\pm\infty} \frac{Q(iy)\re^{-iay-b}}{\sin(iy)}
         = -\frac1{z_0}\mathrm{V.p.}\!\prod_{k=-\infty}^\infty\!\!'
         \Bigl( 1 - \frac{\tilde z_k}{z_k}\Bigr) \ne 0.
\]
On the other hand, the integral representation of $Q$ implies the
estimate $|Q(z)|=\mathrm{o}\bigl(\re^{|\Im z|}\bigr)$ as
$|z|\to\infty$ by~\cite[Lemma~1.3.1]{Ma}, which is inconsistent
with at least one of the above limits. The contradiction derived
shows that $q=0$.
\end{proof}

We denote by $B$ the mapping
\begin{equation}\label{eq:4.1a}
    B\, :\,g \mapsto - A_g^{-1}(\cF^{-1}\sin \be(g)),
\end{equation}
which, according to the above corollary, is well defined in some
neighbourhood of the origin in $X$.

\begin{proof}[Proof of Theorem~\ref{thm:inv}]
Let $g$ be an arbitrary function in $X$ and put $z_n:=\pi n +
e_n(g)$, $n\in\bZ$. Since $e_n(g) \to 0$ as $|n|\to\infty$ by the
Riemann--Lebesgue lemma, there is a number $m_0\in\bN$ such that
$|e_n(g)|<\pi/2$ if $|n|>m_0$. Without loss of generality we may
assume that the sets $\{z_{-m_0},\dots,z_{m_0}\}$ and
$\{z_n\}_{|n|>m_0}$ are disjoint, so that the numbers~$z_n$ with
$|n|>m_0$ occur only once in the sequence $(z_n)_{n\in\bZ}$.

As the system of the functions $\{\re^{2\pi nit}\}_{n\in\bZ}$ is
complete in~$X$ and in $L_1(0,1)$, there exists a natural $m\ge
m_0$ and a trigonometric polynomial $p(t):=\sum_{k=-m}^m c_k
\re^{2\pi kit}$ such that the function $\tilde g:=g-p$ satisfies
the inequalities $\|\tilde g\|_X < (2\rho\|M\|)^{-1}$ and
$\|\tilde g\|_{L_1} < \pi/2$. We set $\zeta _n:=\pi n+ e_n(\tilde
g)$, $n\in\bZ$; then the numbers $\zeta _n$ are pairwise distinct
and $\zeta _n=z_n$ if $|n|>m$. If we put now $\tilde f:= B(\tilde
g)\in X$ with the mapping $B$ defined in~\eqref{eq:4.1a}, then the
function~$G:=F_{\tilde f}$ vanishes at the points $\zeta _n$---in
particular, $G(z_n)=0$ if $|n|>m$.

We take next a sufficiently small number $\eps>0$ such that the
functions $\tilde g_l(t):=\tilde g + \eps \re^{2\pi lit}$,
$l=-m,\dots,m$, satisfy the same norm estimates as $\tilde g$
does. We put $\tilde f_l:=B(\tilde g_l)\in X$ for $l=-m,\dots,m$;
then the corresponding entire functions $G_l:=F_{\tilde f_l}$ have
the following properties:
\begin{equation}\label{eq:4.1}
    G_l(\zeta _l)\ne0
    \quad\text{and}\quad
    G_l(\zeta_k)=0
    \quad\text{if}\quad k\ne l.
\end{equation}

By construction, the functions $G$ and $G_l$, $l=-m,\dots,m$, vanish at
the points $z_n$ for $|n|>m$. We shall show that there exist complex
numbers $\al_{-m}, \dots , \al_m$ such that the function
\[
    F:= G + \sum_{l=-m}^m \al_l (G_l-G)
\]
has also zeros at the points $z_{-m},\dots,z_m$ (of respective
multiplicity). This will prove the theorem, since, as is easily
seen, $F$ is of the form~\eqref{eq:F2t} for the function
\begin{equation}\label{eq:3.f}
    f:=\tilde f + \sum_{l=-m}^m \al_l (\tilde f_l-\tilde f).
\end{equation}

Assume that there is no such numbers $\al_l$, $l=-m,\dots,m$.
Denote by $(w_k)_{k=1}^q$ the subsequence of pairwise distinct
numbers among $z_{-m},\dots,z_m$, and set $r_k$ to be the number
of times $w_k$ appears in $(z_l)_{l=-m}^m$. Then the linear system
of order $2m+1$,
\begin{equation}\label{eq:3.syst}
    \sum_{l=-m}^m \al_l (G_l-G)^{(j)}(w_k)=-G^{(j)}(w_k), \qquad
    k=1,\dots,q, \qquad j=0,\dots,r_k-1,
\end{equation}
has no solution, and hence the corresponding homogeneous system
has a nontrivial solution $\al^0_{-m},\dots,\al^0_m$. Put
$h:=\sum_{l=-m}^m \al^0_l (\tilde f_l-\tilde f)$ and
$H(z):=\int_0^1 h(t)\re^{iz(2t-1)}\,dt$; then the function~$H$
equals $\sum_{l=-m}^m\al_l^0(G_l-G)$ and hence verifies the
relations
\[
    H(z_n) = \int_0^1 h(t)\re^{iz_n(2t-1)}\,dt =0
\]
for $|n|>m$ by the construction of the functions $G$ and $G_l$,
and
\[
    H^{(j)}(w_k)
    =i^j\int_0^1 h(t)(2t-1)^j\re^{iw_k(2t-1)}\,dt =0
\]
for $k=1,\dots,q$ and $j=0,\dots,r_k-1$ by assumption. Since $h\in
X\subset L_1(0,1)$, Lemma~\ref{lem:3.3} yields $h=0$ as in the
proof of Corollary~\ref{cor:4.2} and thus $H\equiv0$. On the other
hand, $H(\zeta_l)= \alpha_l^0 G_l(\zeta_l)$, which is not zero for
at least one $l=-m,\dots,m$ in view of~\eqref{eq:4.1}. The
contradiction derived justifies existence of the required $f\in
X$.

Uniqueness of~$f$ follows by the similar reasoning based on
Lemma~\ref{lem:3.3}.
\end{proof}

\begin{remark}
It can be shown that the mapping $g\mapsto f$ constructed in
Theorem~\ref{thm:inv} is continuous. Indeed, fix a function
$g_0\in X$, take the polynomial $p_0$ as in the proof of that
theorem, and choose a neighbourhood~$\mathcal{O}$ of~$g_0$ such
that the functions $\tilde g := g - p_0$ satisfy $\|\tilde
g\|_X<(2\rho\|M\|)^{-1}$ for all $g\in\mathcal{O}$. Then the
functions $\tilde f$ and $\tilde f_l$, $l=-m,\dots,m$, become
continuous in~$g\in\mathcal{O}$, and in view of
formula~\eqref{eq:3.f} it remains to show that the solution
$\al_{-m}, \dots,\al_m$ of system~\eqref{eq:3.syst} depends
continuously on $g\in\mathcal{O}$. The local continuity is evident
if the numbers $-\pi m + e_{-m}(g),\dots, \pi m + e_m(g)$ are
pairwise distinct. If some of these numbers collide as $g$ tends
to $g_0$, continuation arguments show that $\al_l(g)$,
$l=-m,\dots,m$, have limits, which solve system~\eqref{eq:3.syst}
for $g=g_0$. This justifies the claim.
\end{remark}

\medskip
\textbf{Acknowledgements.} {The first author gratefully
acknowledges the financial support of the Alexander von Humboldt
Foundation and thanks the Institute for Applied Mathematics of
Bonn University for the warm hospitality.}

\end{document}